\documentclass[a4paper,12pt]{article}  
\RequirePackage[OT1]{fontenc}  
\usepackage[natbib,noinfoline]{imsart}  
\usepackage[OT1]{fontenc}  
\usepackage[english]{babel}  
\usepackage[utf8]{inputenc}   
\usepackage{latexsym,amssymb,amsmath,epsfig,psfrag}  
\def\PP{\mathbb{P}} 
\def\RR{\mathbb{R}} 
\def\EE{\mathbb{E}} 
\def\NN{\mathbb{N}} 
\def\ZZ{\mathbb{Z}} 
\def\SS{\mathbb{S}}
 
\def\LL{\mathbf{L}}
\def\Var{\mathsf{Var}}
\def\II{\mbox{ 1\hskip -.29em I}}
\newcommand{\sfrac}[2]{\kern.1em
        \raise.5ex\hbox{$#1$}\kern-.1em
        /\kern-.15em\lower.25ex\hbox{$#2$}}
\def\bdes{\begin{description}}
\def\edes{\end{description}}
\def\ita{\item[(a)]}
\def\itb{\item[(b)]}

\def\iti{\item[(i)]}
\def\itii{\item[(ii)]}
\def\itiii{\item[(iii)]}
\def\itiv{\item[(iv)]}
\def\itv{\item[(v)]}

\newtheorem{defi}{Definition}[section]
\newtheorem{lemm}[defi]{Lemma}
\newtheorem{prop}[defi]{Proposition}
\newtheorem{coro}[defi]{Corollary}
\newtheorem{theo}[defi]{Theorem}
\newtheorem{exem}{Example}
\newtheorem{rem}{Remark}
\newenvironment{dem}{\vskip 2mm\noindent {\it Proof} :}
                    {\hfill $\square$ \vskip 2mm \noindent} 
 
\newcommand{\eps}{\varepsilon}
\def\bex{\begin{exem} \em }
\def\eex{\end{exem} }
\def\brem{\begin{rem} \em}
\def\erem{\end{rem} }

\def\mylabel#1{\label{#1}}
\begin{document}
\begin{frontmatter}

\title{A modified Poincar\'e inequality and its application to First Passage Percolation} 

\runtitle{Modified Poincaré inequality \& FPP}
\thanks{We acknowledge financial support from the Swiss National Science Foundation grant 200021-1036251/1.} 
\author{\fnms{Michel} \snm{Benaïm} \ead[label=e2]{michel.benaim@unine.ch}}
\address{Institut de Mathématiques\\ Universit\'e de Neuchâtel,\\ 11 rue Emile
  Argand,\\ 2000 Neuchâtel, SUISSE\\
\printead{e2}}
\author{\fnms{Rapha\"el} \snm{Rossignol} \ead[label=e1]{raphael.rossignol@unine.ch}}
\address{Institut de Mathématiques\\ Universit\'e de Neuchâtel,\\ 11 rue Emile
  Argand,\\ 2000 Neuchâtel, SUISSE\\
\printead{e1}}

%\affiliation{Université de Neuchâtel}

\runauthor{M. Benaïm, R. Rossignol}

\begin{abstract}:~  We extend a functional inequality for the Gaussian measure on
  $\RR^{n}$ to the one on
  $\RR^{\NN}$. This inequality improves on the classical Poincaré
  inequality for Gaussian measures. As an application, we prove that First Passage Percolation has sublinear
  variance when the edge times distribution belongs to a wide class of
  continuous distributions, including the exponential one. This extends a
  result by Benjamini, Kalai and Schramm \cite{BenjaminiKalaiSchramm03}, valid
  for positive Bernoulli edge times.\vspace*{3mm}

\noindent {\bf Résumé}:~  On étend une inégalité fonctionnelle pour
  la mesure gaussienne sur $\RR^{n}$ à celle sur
  $\RR^{\NN}$. Cette inégalité peut être vue comme une amélioration de l'inégalité de
  Poincaré classique pour la mesure gaussienne. Comme application, nous
  montrons que la percolation de premier passage a une variance sous-linéaire
  pour une large classe de distributions des temps d'arêtes, incluant les
  lois exponentielles. Ceci étend un résultat de Benjamini, Kalai et
  Schramm \cite{BenjaminiKalaiSchramm03}, valable pour des temps d'arêtes strictement positifs suivant une loi de Bernoulli.
\end{abstract}

\begin{keyword}[class=AMS]
\kwd[Primary ]{60E15}
\kwd[; secondary ]{60K35}
\end{keyword}

\begin{keyword}
\kwd{Ornstein-Uehlenbeck process}
\kwd{Poincaré inequality}
\kwd{Hypercontractivity}
\kwd{First Passage Percolation}
\end{keyword}

\end{frontmatter}

\section{Introduction}
\mylabel{sec:intro}
First Passage Percolation was introduced by Hammersley and Welsh
\cite{HammersleyWelsh65} to model the
flow of a fluid in a randomly porous material (see \cite{Howard04} for a recent
account on the subject). We will consider
the following model of First Passage Percolation in $\ZZ^d$, where $d\geq 2$ is an
integer. Let $E=E\left(\ZZ^d\right)$ denote the set of edges in $\ZZ^d$. The
passage time of the fluid through the edge $e$ is denoted by
$x_e$ and is supposed to be nonnegative. Randomness of the porosity is given
by a product probability measure on $\RR_+^E$. Thus,
$\RR_+^E$ is equipped with the probability measure $\mu=\nu^{\otimes E}$,
where $\nu$ is a probability measure on $\RR_+$ according to which each
passage time is distributed, independently from the others. If $u,v$ are two vertices of
$\ZZ^d$, the notation $\alpha:\{u,v\}$ means that $\alpha$ is a path with
end points $u$ and $v$. When $x\in\RR_+^E$, $d_x(u,v)$ denotes the first passage
time, or equivalently the distance from $u$ to $v$ in the metric induced by $x$,
$$d_x(u,v)=\inf_{\alpha:\{u,v\}}\sum_{e\in\alpha}x_e\;.$$
The study of $d_x(0,nu)$ when $n$ is an integer
which goes to infinity is of central importance. Kingman's subadditive
ergodic theorem implies the existence, for each fixed $u$, of a ``time
constant'' $t(u)$ such that:
$$\frac{d_x(0,nu)}{n}\xrightarrow[n\rightarrow +\infty]{\nu-a.s} t(u)\;.$$
It is known (see Kesten \cite{Kesten84}, p.127 and 129) that if $\nu(\{0\})$ is
strictly smaller than the critical probability for Bernoulli bond percolation
on $\ZZ^d$, then $t(u)$ is  positive for every $u$ which is not the
origin. Thus, in this case, one can say that the random variable $d_x(0,nu)$
is located around $nt(u)$, which is of order $O(|nu|)$, where we denote by $|.|$ the $L^1$-norm of vertices in $\ZZ^d$. In this paper, we
are interested in the fluctuations of this quantity. Precisely, we define, for
any vertex $v$,
$$\forall x\in\RR_+^E,\; f_v(x)=d_x(0,v)\;.$$
It is widely believed that the fluctuations of $f_v$ are of order
$|v|^{1/3}$. Apart from some predictions made by physicists, this faith relies on recent results for related growth models \cite{BaikDeiftJohansson99,Johansson00,Johansson00phys}. Until recently, the best
results rigourously obtained for the fluctuations of $f_v$ were some
moderate deviation estimates of order $O(|v|^{1/2})$ (see \cite{Kesten93,Talagrand95}). In 2003, Benjamini, Kalai and Schramm
\cite{BenjaminiKalaiSchramm03} proved that for Bernoulli
edge times, the variance of $f_v$ is of order $O(|v|/\log|v|)$, and therefore,
the fluctuations are of order  $O(|v|^{1/2}/(\log|v|)^{1/2})$.

The main result of
the present paper is Theorem \ref{theoFPPexpo}, where we extend the result
of Benjamini et al. to a large class of probability measures,
absolutely continuous with respect to the Lebesgue measure. This class includes
all the gamma and beta
distributions (and therefore, the exponential distribution), but also any probability measure whose
density is bounded away from 0 on its support, and notably the uniform
distribution on
$[a,b]$, with $0\leq a<b$. The result of Benjamini et al. relies mainly on an
inequality on the discrete cube due to Talagrand \cite{Talagrand94a}. In their
article, they suggested to extend their work to other edge-times distributions
by using the tools developped by Ledoux \cite{Ledoux01}, and they pointed out a Gaussian
version of Talagrand's inequality found by Bobkov and Houdr\'e
\cite{BobkovHoudre99}. Our strategy to extend Benjamini, Kalai and Schramm's
result is thus to start from the Gaussian version of Talagrand's inequality, and then to adapt their argument to a continuous
distribution $\nu$ via this inequality and a standard change of variable from
the Gaussian distribution to $\nu$.  Following Ledoux \cite{LedouxJerusalem}, p. 41,
we shall give another proof of
the Gaussian inequality that we need in Proposition
\ref{propTaladaptFPPN}. Since it follows (almost) entirely Ledoux's argument, and since
the result is implied, up to a multiplicative constant, by Bobkov and Houdr\'e
\cite{BobkovHoudre99}, we don't claim that this proposition is original at
all. Nevertheless, we choosed to write it down, because in the precise form which we
state it in, it really improves upon the classical Poincaré inequality for the
Gaussian measure. By this, we mean that it implies the classical Poincaré
inequality with the true, optimal constant.

This article is organized as follows. In Section \ref{sec:talagrand}, we derive the Gaussian
analog to Talagrand's inequality \cite{Talagrand94a}, and extend it to a
countable product of Gaussian measures. Section \ref{sec:extension} is devoted to the obtention
of similar inequalities for other continuous measures by a simple mean of
change of variable. This allows us to adapt the argument of Benjamini et
al. to some continuous settings in section \ref{sec:FPP}.

\paragraph{Notation}
Given a probability space $(\Omega, {\cal A}, \mu)$ and a real valued measurable function $f$ defined on $\Omega$ we let  $$\|f\|_{p,\mu} = \left(\int \|f\|^p\;d\mu\right)^{1/p} \in [0, \infty]$$ and we let $L^p(\mu)$ denote the set of $f$ such that $\|f\|_{p,\mu} < \infty.$ 
The {\em mean} of $f \in L^1(\mu)$ is denoted
$$\EE_{\mu}(f)  = \int f d\mu$$ and  the {\em variance} of $f \in L^2(\mu)$ is $$\Var_{\mu}(f) = \|f - \EE_{\mu}(f)\|_{2,\mu}^2.$$
When the choice of $\mu$ is unambiguous we may write $\|f\|_p$ (respectively  $L^p$, $\EE(f)$, $\Var(f)$)    for $\|f\|_{p,\mu}$  (respectively  $L^p(\mu)$, $\EE_{\mu}(f)$, $\Var_{\mu}(f)$). 

\section{A modified Poincar\'e inequality for Gaussian measures on $\RR^{\NN}$}
\mylabel{sec:talagrand}
The purpose of this section is to rewrite Ledoux's argument (see \cite{LedouxJerusalem}) in deriving a Poincar\'e type inequality for a
 product measure having the form $\lambda \otimes \gamma^{\NN}$ where
 $\lambda$ is uniform on $\{0,1\}^S$ with $S$ finite and $\gamma^{\NN}$ is the
 standard  Gaussian measure on $\RR^{\NN}.$ Such an inequality, which will prove to be very useful in the sequel, can be seen as
 a ``Gaussian'' version of an inequality proved by Talagrand in
 \cite{Talagrand94a} for Bernoulli measures. Please note that it is implied,
 up to a multiplicative constant, by an inequality from Bobkov and Houdr\'e \cite{BobkovHoudre99}. We shall first prove it for the standard Gaussian on $\RR^n$ and then extend it to the case of $\RR^{\NN}.$
\subsection*{{\em The case of $\RR^n$}}
We let  
$$\gamma(dy) = \frac{1}{\sqrt{2\pi}} e^{-y^2/2} dy$$
denote the  standard Gaussian measure on $\RR$  
and $$\beta = \frac{1}{2}(\delta_0 + \delta_1)$$ the Bernoulli measure on $\{0,1\}.$ If $n \in \NN,$ and $S$ is a finite set we let $\gamma^n = \gamma^{\otimes^n}$ and  $\lambda = \beta^{\otimes^S}$  denote the associated product measures  on $\RR^n$ and $\{0,1\}^S.$ 

Given a  measurable mapping $f :  \{0,1\}^S\times\RR^n \mapsto \RR$  and $q \in S$ we set
 $$\nabla_q f(x,y)=f(x,y)-\int_{\{0,1\}} f(x,y)d\lambda (x_q).$$
A map $f$ is said to be  {\em weakly differentiable} (in the $y$ variables) provided there exist locally integrable (in the $y$ variables) functions denoted  $(\frac{\partial f}{\partial y_j}),  j = 1, \ldots, n,$ such that
$$\int_{\RR^{n}} \frac{\partial f}{\partial y_j}(x,y) g(y) dy = - \int_{\RR^{n}} \frac{\partial g}{\partial y_j}(y) f(x,y) dy$$
for every smooth function $g : \RR^n \mapsto \RR$ with compact support.

The {\em weighted Sobolev} space $H_1^2\left(\lambda\otimes \gamma^n\right)$ is defined to be the space of weakly differentiable functions $f$ on $\{0,1\}^S\times\RR^n$ 
such that
$$\|f\|^2_{H_1^2} = \|f\|_2^2 + \sum_{j =1}^n \left\|\frac{\partial f}{\partial
    y_j} \right\|_2^2 < \infty.$$

The
following result is essentially a Gaussian version of Talagrand's Theorem 1.5 in \cite{Talagrand94a}.
\begin{prop}
\mylabel{propTaladaptFPPn}
For every function $f$ in $H_1^2(\lambda\otimes\gamma^{n})$,
$$\Var(f)\leq \sum_{q\in S}\left\|\nabla_q f\right\|_2^2+\sum_{j=1}^n\left\|\frac{\partial f}{\partial
    y_i}\right\|_2^2 \phi\left(\frac{\left\|\frac{\partial f}{\partial
    y_i}\right\|_1}{\left\|\frac{\partial f}{\partial
    y_i}\right\|_2}\right)\;,
$$
where $\phi : [0,1] \mapsto [0,1]$ is defined as  $$\phi(u) = 2 \int_0^1 \frac{u^{2t}}{(1+t)^2} dt.$$  
\end{prop}
\brem
 Function $\phi$ is continuous nondecreasing   with $\phi(0) = 0$ (in fact 
$\phi(u) \sim_{u = 0} -1/(\log(u))$ ) and $\phi(1) = 1.$ In particular, Proposition \ref{propTaladaptFPPn} implies the standard Poincar\'e inequality for $\gamma^n.$ That is
$$\Var(f)\leq \sum_{i=1}^n\left\|\frac{\partial f}{\partial
    y_i}\right\|_2^2\;,$$
for all $f \in H^1_2(\gamma^n).$
\erem
\begin{dem}
 For every  $f \in L^2$ 
\begin{equation} 
\label{eq2var} \Var_{\lambda\otimes\gamma^{n}}(f)=\EE_{\gamma^n}(\Var_{\lambda}(f))+\Var_{\gamma^n}(\EE_{\lambda}(f)).
\end{equation}
Furthermore, 
variance satisfies the following well known tensorisation property (see for
instance Ledoux \cite{Ledoux01}, Proposition 5.6 p. 98):
$$\forall g\in L^2\left(\lambda\right),\;\Var_{\lambda}(g)\leq \sum_{q\in
  S}\left\|\nabla_q g\right\|_{2,\lambda}^2\;.$$
Hence 
\begin{equation}
\label{eqmajovarlambda}
\EE_{\gamma^n}(\Var_{\lambda}(f))\leq \sum_{q\in
  S}\left\|\nabla_q f\right\|_{2,\lambda\otimes\gamma^{n}}^2\;.
\end{equation}
Following Benjamini et al. \cite{BenjaminiKalaiSchramm03}, we shall use an 
hypercontractivity  property of a semi-group $(P_t)_{t\geq
  0}$  with invariant measure $\gamma^n$ in order to bound the last term in
the right hand side of inequality
(\ref{eq2var}). A most natural choice is to take for  $(P_t)_{t\geq
  0}$ the semi-group induced by 
$n$ independent copies of an
Ornstein-Uehlenbeck process.

Let $C_0^{\infty}(\RR^n)$ be  the space of smooth real valued functions on  $\RR^n$ that go to zero at infinity and let $\mathcal{A} \subset C_0^{\infty}(\RR^n)$ be the set of functions whose partial derivatives of all order $\frac{\partial^k f}{\partial x_{i_1} \ldots x_{i_k}}$ lie in $C_0^{\infty}(\RR^n).$ 
For any $g \in \mathcal{A}$ and $t \geq 0$,
$$\frac{dP_t(g)}{dt} =  \LL P_t(g)\;,$$
where 
$$\LL g(x)=\Delta g(x)-x.\nabla g(x)\;.$$
 If $f$ is defined on $\{0,1\}^S\times\RR^n$ and $x\in\{0,1\}^S$, we let  $f_x$ denote the function on $\RR^n$ defined as 
$$\forall y\in \RR^n,\;f_x(y)=f(x,y)\;.$$
Let $\mathcal{A}_S$ be the set of functions $f$ on $\{0,1\}^S\times\RR^n$ such
that, for every $x\in\{0,1\}^S$, the function $f_x$ belongs to $\mathcal{A}$.\\
Without loss of generality we may (hence do) assume that $f \in 
\mathcal{A}_S.$ For, by standard approximation results  $\mathcal{A}_S$ is a dense subset of $H_1^2(\lambda\otimes\gamma^{n})$ equipped with the norm $\|f\|_{H_1^2}$  and  the inequality to be proved is a closed condition in  $H_1^2(\lambda\otimes\gamma^{n}).$ 

Given such an $f$, let $F = \EE_{\lambda}(f)$, which belongs to ${\cal A}.$ Then
\begin{eqnarray*}
P_0F&=&F\;,\\
P_tF&\xrightarrow[t\rightarrow \infty]{L^2}& \int F(y)\;d\gamma^n(y)\;.
\end{eqnarray*}
Therefore,
\begin{eqnarray*}
\Var_{\gamma^n}\left(F\right)&=&\int_0^{+\infty}-\frac{d}{dt}\EE_{\gamma^n}((P_tF)^2)\;dt,\;\\
&=& 2\int_0^{+\infty}\EE_{\gamma^n}(-P_tF\LL P_tF)\;dt.\;\\
\end{eqnarray*}
Notice that for every function $g$ in $\mathcal{A}$, it follows from
integration by parts that:
\begin{eqnarray*}
 \EE_{\gamma^n}(-g\LL g)&=&\sum_{i=1}^n\left\|\frac{\partial g}{\partial
      y_i}\right\|_2^2\;.
\end{eqnarray*}
Consequently,
\begin{equation}
\label{eqvarF}
\Var_{\gamma^n}(F)=2\sum_{i=1}^n\int_0^{+\infty}\EE_{\gamma^n}\left\lbrack\left(\frac{\partial}{\partial
      y_i}P_tF\right)^2\right\rbrack\;dt.\;
\end{equation}
Recall that, by Mehler's formula
$$P_tg(y)=\int g\left(ye^{-t}+z\sqrt{1-e^{-2t}}\right)\;d\gamma_n(z)\;.$$
This expression implies that, for any $i=1,\ldots ,n$ and $g \in {\cal A}$,
$$\frac{\partial}{\partial y_i}P_tg=e^{-t}P_t\left(\frac{\partial g}{\partial
    y_i}\right)\;.$$
Thus,
$$\EE_{\gamma^n}\left\lbrack\left(\frac{\partial}{\partial
      y_i}P_tF\right)^2\right\rbrack=e^{-2t}\EE_{\gamma^n}\left\lbrack\left(P_t\frac{\partial
      F}{\partial
      y_i}\right)^2\right\rbrack\;.$$
Nelson's Theorem asserts that $(P_t)_{t\geq
  0}$ is a hypercontractive semi-group with hypercontractive function
  $q(t)=1+\exp(2t)$ (see Nelson \cite{Nelson73} or An\'e et al. 
  \cite{thesardsLedoux} p.~22):
$$\forall t\geq 0, \; \forall g\in
L^2\left(\gamma^n\right),\;\|P_tg\|_{2}\leq \|g\|_{q^*(t)}\;,$$
where $q^*(t)=1+\exp(-2t)$ is the conjugate exponent of $q(t)$. Thus,
\begin{eqnarray}
\nonumber \EE_{\gamma^n}\left\lbrack\left(\frac{\partial}{\partial
      y_i}P_tF\right)^2\right\rbrack&\leq
      &e^{-2t}\left\lbrack\EE_{\gamma^n}\left(\frac{\partial F}{\partial
      y_i} \right)^{q^*(t)}\right\rbrack^{\frac{2}{q^*(t)}}\;,\\
\nonumber &=&e^{-2t}\left\lbrack\EE_{\gamma^n}\left(\int\frac{\partial f}{\partial
      y_i}\;d\lambda\right)^{q^*(t)}\right\rbrack^{\frac{2}{q^*(t)}}\;,\\
\label{majnelson} &\leq&e^{-2t}\left\lbrack\EE_{\lambda\otimes\gamma^n}\left(\int\frac{\partial
      f}{\partial
      y_i}\right)^{q^*(t)}\right\rbrack^{\frac{2}{q^*(t)}}\;,
\end{eqnarray}
which follows  from Jensen's inequality. Now, according to  H\"older inequality,
$$\EE\left\lbrack\left(\frac{\partial f}{\partial
    y_i}\right)^{q^*(t)}\right\rbrack\leq \EE\left\lbrack\left(\frac{\partial f}{\partial
    y_i}\right)^{2}\right\rbrack^{q^*(t)-1}\EE\left\lbrack\left|\frac{\partial f}{\partial
    y_i}\right|\right\rbrack^{2-q^*(t)}.\;$$
Thus, inequality (\ref{majnelson}) implies:
\begin{eqnarray*}
&&\int_0^{+\infty}\EE_{\gamma^n}\left\lbrack\left(\frac{\partial}{\partial
      y_i}P_tF\right)^2\right\rbrack\;dt\\
&\leq &\int_0^{+\infty}e^{-2t}\left\lbrack\EE\left\lbrack\left(\frac{\partial f}{\partial
    y_i}\right)^{2}\right\rbrack^{q^*(t)-1}\EE\left\lbrack\left|\frac{\partial f}{\partial
    y_i}\right|\right\rbrack^{2-q^*(t)}\right\rbrack^{\frac{2}{q^*(t)}}\;dt\;,\\
&=& \left\|\frac{\partial f}{\partial
    y_i}\right\|_2^2\int_0^{+\infty}e^{-2t}\left(\frac{\left\|\frac{\partial f}{\partial
    y_i}\right\|_1}{\left\|\frac{\partial f}{\partial
    y_i}\right\|_2}\right)^{\frac{2(2-q^*(t))}{q^*(t)}}\;dt\;,\\
&=&
\left\|\frac{\partial f}{\partial
    y_i}\right\|_2^2\int_0^{1}\left(\frac{\left\|\frac{\partial f}{\partial
    y_i}\right\|_1}{\left\|\frac{\partial f}{\partial
    y_i}\right\|_2}\right)^{2s}\frac{1}{(1+s)^2}\;ds.\
\end{eqnarray*}
where last equality follows from the change of variable 
$s = \frac{1-e^{-2t}}{1+e^{-2t}}.$
A combination of this upper bound with equality (\ref{eqvarF}) and inequalities (\ref{eqmajovarlambda}) and 
 (\ref{eq2var}) gives the desired result.
\end{dem}
\subsection*{{\em The case of $\RR^{\NN}$}}
Let now $H_1^2\left(\lambda\otimes \gamma^\NN\right)$ be the space of functions $f \in L^2\left(\lambda\otimes \gamma^\NN\right)$ verifying the two following conditions:
\begin{description}
\ita  For all $i\in \NN,$
there exists a function $h_i$ in $L^2\left(\lambda\otimes \gamma^\NN\right)$
such that $$-\int_{\RR} g'(y_i)f(x,y)\;dy_i=\int_{\RR} g(y_i)h_i(x,y)\;dy_i\;, \lambda\otimes \gamma^\NN \; a.s $$
for every smooth function $g :\RR \mapsto \RR$ having compact support.
The function $h_i$ is called the partial derivative of $f$ with respect to $y_i$,
and is denoted by $\frac{\partial f}{\partial
    y_i}.$
\itb  The  sum of the $L^2$ norms of the partial derivatives of $f$ is finite:
$$\sum_{i\in\NN}\left\|\frac{\partial f}{\partial
    y_i}\right\|_2^2<\infty\;.$$
\end{description}
\begin{rem}
It is not hard to verify that for  $f$ depending on finitely many variables, say $y_1, \ldots, y_n,$ then  $f \in H_1^2(\lambda\otimes\gamma^{\NN}) \Leftrightarrow f \in H_1^2(\lambda\otimes\gamma^{n}).$
\end{rem}
Let this Sobolev space be equipped with the norm:
$$\|f\|_{H}=\sqrt{\|f\|_2^2+\sum_{q\in S}\|\nabla_qf\|_2^2+
\sum_{i\in\NN}\left\| \frac{\partial f}{\partial
    y_i}\right\|_2^2}\;.$$
Proposition \ref{propTaladaptFPPn} extends to functions of countably infinite
Gaussian variables as follows:
\begin{prop}
\mylabel{propTaladaptFPPN}
For every function $f$ in $H_1^2(\lambda\otimes\gamma^{\NN})$,
$$\Var(f)\leq \sum_{q\in S}\left\|\nabla_q f\right\|_2^2+\sum_{i\in\NN}\left\|\frac{\partial f}{\partial
    y_i}\right\|_2^2 \phi\left(\frac{\left\|\frac{\partial f}{\partial
    y_i}\right\|_1}{\left\|\frac{\partial f}{\partial
    y_i}\right\|_2}\right)\;,$$
where $\phi$ is as in Proposition \ref{propTaladaptFPPn}.\\
\end{prop}
The proof of Proposition \ref{propTaladaptFPPN} relies on  
the following simple approximation lemma combined with  Proposition \ref{propTaladaptFPPn}.
\begin{lemm}
\label{lemconvH12}
 Let $\mathcal{F}_n$ be the $\sigma$-algebra generated by the first $n$
coordinate functions in $\RR^\NN$. Let $f \in H_1^2\left(\lambda\otimes
  \gamma^\NN\right)$ and $f_n = \EE\left(f|\mathcal{F}_n\right)$ be the
conditional expectation of $f$ with respect to $\mathcal{F}_n$. Then
\begin{enumerate}
\iti $f_n \in H_1^2\left(\lambda\otimes \gamma^n\right)$,
\itii For every $i$ in $\{1,\ldots ,n\}$,
$$\frac{\partial f_n}{\partial y_i}=\EE\left(\frac{\partial f}{\partial y_i}|\mathcal{F}_n\right)\;,$$
\itiii The following convergence takes place in $H_1^2\left(\lambda\otimes
    \gamma^\NN\right)$:
$$f_n \xrightarrow[n\rightarrow \infty]{}f\;.$$
In particular,
\begin{equation}
\label{eqlim1}
\sum_{q\in S}\left\|\nabla_q f_n\right\|_2^2\xrightarrow[n\rightarrow
+\infty]{}\sum_{q\in S}\left\|\nabla_q f\right\|_2^2,
\end{equation}
\begin{equation}
\label{eqlim2}
\Var(f_n)\xrightarrow[n\rightarrow
+\infty]{}\Var(f),
\end{equation}
and
\begin{equation}
\label{eqlim3}
\sum_{i=1}^n\left\|\frac{\partial f_n}{\partial
    y_i}\right\|_2^2 \phi\left(\frac{\left\|\frac{\partial f_n}{\partial
    y_i}\right\|_1}{\left\|\frac{\partial f_n}{\partial
    y_i}\right\|_2}\right) \xrightarrow[n\rightarrow +\infty]{}\sum_{i\in\NN}\left\|\frac{\partial f}{\partial
    y_i}\right\|_2^2 \phi\left(\frac{\left\|\frac{\partial f}{\partial
    y_i}\right\|_1}{\left\|\frac{\partial f}{\partial
    y_i}\right\|_2}\right)\;.
\end{equation}
\end{enumerate}
\end{lemm}
\begin{dem}
Of course, $\EE\left(f|\mathcal{F}_n\right)$ belongs
  to $L^2\left(\lambda\otimes \gamma^n\right)$. Let
 $g : \RR^n \mapsto \RR$ be a smooth function with compact support 
  and $i \in \{1,\ldots,n\}$. According to Fubini's Theorem,
\begin{eqnarray*}
&& -\int\frac{\partial g}{\partial
 y_i}\EE\left(f|\mathcal{F}_n\right)\;dy_1\ldots dy_n\\
 &=& -\int\EE\left(\int \frac{\partial g}{\partial
 y_i}f\;dy_i|\mathcal{F}_n\right)\;dy_1\ldots dy_{i-1}dy_{i+1} \ldots dy_n
 \;.
\end{eqnarray*}
Then, it follows from the definition of $H_1^2\left(\lambda\otimes
 \gamma^\NN\right)$ that:
$$-\int \frac{\partial g}{\partial y_i}f\;dy_i=\int g \frac{\partial
  f}{\partial y_i}\;dy_i\;,$$
with probability 1. Therefore,
$$-\int\frac{\partial g}{\partial
 y_i}\EE\left(f|\mathcal{F}_n\right)\;dy_1\ldots dy_n =\int g
 \EE\left(\frac{\partial f}{\partial
 y_i}|\mathcal{F}_n\right)\;dy_1\ldots dy_n\;.$$
proving assertions $(i)$ and $(ii)$ of the Lemma.

Since for any function $g$ in $L^2\left(\lambda\otimes
  \gamma^\NN\right)$, the martingale  $(\EE\left(g|\mathcal{F}_n\right))_{n \geq 0}$ converges to $g$ in
$L^2\left(\lambda\otimes \gamma^\NN\right),$ one has
$$\EE\left(\EE\left(\frac{\partial f}{\partial y_i}|\mathcal{F}_n\right)^2\right)\xrightarrow[n\rightarrow \infty]{}\EE\left(\frac{\partial f}{\partial y_i}\right)\;.$$
In addition, by Jensen's inequality,
$$\EE\left(\EE\left(\frac{\partial f}{\partial y_i}|\mathcal{F}_n\right)\right)^2\leq
 \EE\left(\left(\frac{\partial f}{\partial y_i}\right)^2\right)\;,$$
and since $f$ belongs to $H_1^2\left(\lambda\times  \gamma^\NN\right)$,
$$\sum_{i\in\NN}\left\|\frac{\partial f}{\partial
    y_i}\right\|_2^2 <\infty\;.$$
By Lebesgue convergence theorem, it follows that:
\begin{equation}
\label{eqdomcon}
\sum_{i\in\NN}\left\|\frac{\partial \EE\left(f|\mathcal{F}_n\right)}{\partial
    y_i}\right\|_2^2\xrightarrow[n\rightarrow \infty]{}\sum_{i\in\NN}\left\|\frac{\partial f}{\partial
    y_i}\right\|_2^2 <\infty\;.
\end{equation}
Thus,
$$\sum_{i\in\NN}\left\|\frac{\partial }{\partial
    y_i}\left\lbrack f-\EE\left(f|\mathcal{F}_n\right)\right\rbrack\right\|_2^2\xrightarrow[n\rightarrow
    \infty]{}0\;.$$
Finally, it is trivial to chek that $\nabla_q\EE\left(f|\mathcal{F}_n\right)$
    converges to $\nabla_qf$ in $L^2\left(\lambda\times \gamma^\NN\right)$ as $n$
    tends to infinity. Therefore, the convergence takes place in $H_1^2\left(\lambda\times  \gamma^\NN\right)$.
Assertions (\ref{eqlim1}) and (\ref{eqlim2}) follow while the  proof of
    assertion (\ref{eqlim3}) follows from Lebesgue convergence theorem just as the proof
 of (\ref{eqdomcon}).
\end{dem}

\section{Extension to other measures}
\mylabel{sec:extension}
As usual, we can deduce from Proposition \ref{propTaladaptFPPN} other
inequalities by mean of change of variables. To make this precise, let $\Omega$ be a measurable space and $\Psi : \RR^{\NN} \mapsto \Omega$ a measurable isomorphism  (meaning that $\Psi$ is one to one with $\Psi$ and $\Psi^{-1}$ measurables). Let $\Psi^*\gamma^{\NN}$ denote the image of $\gamma^{\NN}$ by $\Psi.$ That is $\Psi^*\gamma^{\NN}(A) = \gamma^{\NN}(\Psi^{-1}(A)).$ For $g : S \times \Omega \mapsto \RR$ such that $g \circ (Id,\Psi) \in H_1^2(\lambda \otimes \gamma^{\NN}),$ one obviously has
 $$\Var_{\lambda \otimes \Psi^*\gamma^{\NN}}(g) = \Var_{\lambda \otimes \gamma^{\NN}}(g \circ \Psi)$$
and
$$\|\partial_{i,\Psi} g\|_{p,\Psi^*\gamma^{\NN}} = \left\|\frac{\partial g \circ \Psi}{\partial y_i}\right\|_{p,\gamma^{\NN}}$$
where $\partial_{i,\Psi} g$ is defined as
$$\partial_{i,\Psi} g(x,\omega) = \frac{\partial (x \circ \Psi)}{\partial y_i}(q,\Psi^{-1}(\omega)).$$ Hence inequality in Proposition \ref{propTaladaptFPPN}  for $f = g \circ (Id,\Psi)$ transfers to 
the same inequality for $g$ provided
$\frac{\partial f}{\partial y_i}$ is replaced  by $\partial_{i,\Psi} g.$

\bex
\mylabel{ex:unisphere}
Let $k \geq 2$ be an integer, $\SS^{k-1} \subset \RR^k$ the unit $k-1$ dimensional sphere, $\Omega = (\RR^+_{*} \times \SS^{k-1})^{\NN},$ and let  $E = (\RR^k_{*})^{\NN}.$ A typical point in $\Omega$ will be written as $(\rho,\theta) = (\rho^i,\theta^i)$ and a typical point in $E$ as $y = (y^j).$ Now consider the change of variables 
$\Psi : E \mapsto \Omega$  given by $\Psi(y) = (\Psi(y)^j)$  with 
$$\Psi^j(y) = \left(\|y^j\|^2, \frac{y^j}{\|y^j\|}\right).$$ The image of $(\gamma^k)^{\NN} = \gamma^{\NN}$ by $\Psi$ is the product measure $\tilde{\gamma}^{\NN}$
where 
$\tilde{\gamma}$ is the probability measure  on $\RR^+_{*} \times S^{k-1}$ defined by
$$ \tilde{\gamma}(dt dv) =  \frac{1}{r_k}e^{-t/2} t^{k/2 -1}\mathbf{1}_{t > 0} dt dv$$
Here $r_k = \int_0^{\infty}  e^{-t/2} t^{k/2 -1} dt$, and $dv$ stands for the
uniform probability  measure on $\SS^{k-1}$. For $g : \{0,1\}^S \times \Omega
\mapsto \RR$ with $g \circ (Id,\Psi) \in H_1^2(\lambda\otimes\gamma^{\NN})$, let 
$$\left (\frac{\partial g}{\partial \rho^j}(x,\rho, \theta), \nabla_{\theta^j}
  g(x,\rho,\theta) \right ) \in \RR \times T_{\theta^i} \SS^{k-1} \subset \RR
\times \RR^k$$ denote the partial gradient of $g$ with respect to the variable
$(\rho^i, \theta^i)$ where $T_{\theta^i} S^{k-1} \subset \RR^k$ stands for the tangent space of $\SS^{k-1}$ at $\theta_i.$
It is not hard to verify that for all $i \in \NN$ and $j \in \{1, \ldots, k\}$,
\begin{equation}
\label{eq:gradsphere}
\partial_{i,j,\Psi} g (x,\rho, \theta) 
= 2 \frac{\partial g}{\partial \rho^i}(x,\rho, \theta) \sqrt{\rho^i} \theta^i_j 
+ \frac{1}{\sqrt{\rho^i}} [\nabla_{\theta^i} g(x,\rho,\theta)]_j.
\end{equation}
\eex
\begin{coro}[$\chi^2$ distribution]
\label{coroTalFPPexpo}
Let $\alpha > 0,\; k \geq 2$ an integer,  and let $\RR^+_*$ be equipped with
the  distribution 
$$\nu(dt) = \frac{e^{-\alpha t} t^{k/2 -1}}{\int_0^{\infty} e^{-\alpha s}
  s^{k/2 -1} ds}\mathbf{1}_{t > 0} dt\;,$$
  and ${\RR^+_*}^{\NN}$ with the product measure $\nu^{\NN}$.
For $g \in H_1^2(\nu^{\NN})$ 
define 
$${\nabla}_ig(y) =  \frac{\partial g}{\partial y_i}(y) \sqrt{y_i}\;.$$
Then 
\begin{equation}
\label{eq:chi2}
\Var(g) \leq
%\sum_{q\in S}\left\|\nabla_q g\right\|_2^2 +
\frac{2}{\alpha}\sum_{i \in \NN} \|{\nabla}_ig\|_2^2 \phi\left( c(k) \frac{  \|{\nabla}_ig\|_1}{\|{\nabla}_ig\|_2}\right)
\end{equation}
with $$c(k) = \sqrt{k} \frac{ \int_0^{\pi}
  |cos(t)|\sin(t)^{k-2}dt}{\int_0^{\pi} \sin(t)^{k-2}dt} = \frac{2 \sqrt{k}}{(k-1)\int_0^{\pi} \sin(t)^{k-2}dt}.$$
\end{coro}
\begin{dem}
follows from (\ref{eq:gradsphere}) applied to the map $(x,\rho,\theta) \rightarrow  g(\frac{\rho}{2\alpha}).$ Details are left to the reader.
\end{dem}
\begin{coro}[Uniform distribution on $\SS^{n}$]
Let $dv_n$ denote the normalized Riemannian probability measure on $\SS^{n} \subset \RR^{n+1}.$ For $g \in H^1_2(dv_n)$  and $i = 1, \ldots, n+1$ let $\nabla_i g(\theta)$ denote the $i^{th}$ component of $\nabla g(\theta)$ in $\RR^{n+1}$ (we see $T_{\theta} \SS^{n}$ as the vector space of $\RR^{n+1}$ consisting of vector that are orthogonal to $\theta$). Then, for $n \geq 2,$
\begin{equation}
\label{eq:sphere}
\Var(g) \leq
%\sum_{q\in S}\left\|\nabla_q g\right\|_2^2 +
\frac{1}{n-1}\sum_{i \in \NN} \|{\nabla}_ig\|_2^2 \phi\left(\frac{  \|{\nabla}_ig\|_1}{\|{\nabla}_ig\|_2}\right)\;.
\end{equation}
\end{coro}
\begin{dem}
follows from (\ref{eq:gradsphere}) applied to the map $(x,\rho,\theta) \rightarrow  g(\theta).$ Details are left to the reader.
\end{dem}
\bex
\mylabel{ex:changementvar}
If one wants to get a result similar to Proposition \ref{propTaladaptFPPN} with $\gamma$
replaced by another probability measure $\nu$, one may of course perform the
usual change of variables through inverse of repartition function. In the
sequel, we denote by
\begin{equation}
\label{defg}
g(x)=\frac{1}{\sqrt{2\pi}}e^{-\frac{x^2}{2}}
\end{equation}
the density of the normalized Gaussian distribution, and by
\begin{equation}
\label{defG}
G(x)=\int_{-\infty }^xg(u)\;du
\end{equation}
its repartition function. For any function $\phi$ from $\RR$ to $\RR$, we
shall note $\tilde{\phi}$ the function from $\RR^\NN$
to $\RR^\NN$ such that $(\tilde{\phi}(x))_j=\phi(x_j)$.
\begin{coro}[Unidimensional change of variables]
\label{corochgtvar}
 Let $\nu$ be a probability on  $\RR^+$  absolutely continuous with respect to the  Lebesgue
measure, with density $h$ and repartition function
$$H(t)=\int_0^th(u)\;du\;,$$
then, for every function $f$ on $\{0,1\}^S\times\RR^n$ 
such that $f\circ (Id, \widetilde{H^{-1}\circ G})\in H_1^2\left(\lambda\otimes \gamma^\NN\right)$,
$$Var_{\lambda_S\otimes\nu^{\otimes \NN}}(f)\leq\sum_{q\in S}\left\|\nabla_q f\right\|_2^2+ 2
\sum_{i\in \NN}\left\|\nabla_if\right\|_2^2\phi\left(\frac{\left\|\nabla_i
      f\right\|_1}{\left\|\nabla_if\right\|_2}\right)\;,
$$
where for every integer $i$,
$$\nabla_if(x,y)=\psi(y_i)\frac{\partial
  f}{\partial y_i}(x,y)\;,$$
and $\psi$ is defined on $I=\{t\geq 0\mbox{ s.t. }h(t)>0\}$:
$$\forall t\in I,\;\psi(t)=\frac{g\circ G^{-1}(H(t))}{h(t)}\;.$$
\end{coro}
\begin{dem}
It is a straightforward consequence of Proposition \ref{propTaladaptFPPN}, applied to $f\circ (Id, \widetilde{H^{-1}\circ G})$.
\end{dem}
\eex

%\begin{coro}[Exponential distribution]
%\label{coroTalFPPexpo}
%Let $\RR^+_*$ be equipped with the $exponential$ distribution $\nu(dt) = e^{-t}\mathbf{1}_{t > 0} dt$ and ${\RR^+_*}^{\NN}$ with the product measure $\nu^{\NN}$.
%For $g \in H_1^2(\lambda \otimes \nu^{\NN})$ 
%define 
%$$\tilde{\nabla}_ig(y) = \frac{\partial g}{\partial y_i}(y) \sqrt{y_i}.$$
%Then 
%$$
%\Var(g)\leq\sum_{q\in S}\left\|\nabla_q g\right\|_2^2+
%\sum_{i\in\NN}\left\|\nabla_ig\right\|_2^2\phi\left(\frac{2\sqrt{2}\left\|\nabla_i
%      g\right\|_1}{\pi\left\|\nabla_ig\right\|_2}\right)\;.
%$$
%\end{coro}
%\begin{dem} Follows from corollary (\ref{chi2}) applied to the map $x \mapsto g(x/2).$
%\end{dem} 

\section{Application to First Passage Percolation}
\mylabel{sec:FPP}
It turns out that Corollary \ref{coroTalFPPexpo} is particularly well suited
to adapt the argument of Benjamini, Kalai and Schramm
\cite{BenjaminiKalaiSchramm03} to show that First Passage Percolation has
sublinear variance when the edges have a $\chi^2$ distribution. This includes
the important case of exponential distribution, for which First Passage
Percolation becomes equivalent to a version of Eden growth model (see for instance Kesten
\cite{Kesten84}, p.130). We do not want
to restrict ourselves to those distributions. Nevertheless, due to the particular strategy that we
adopt, we can only prove the result for some continuous edge distributions
which behave roughly like a gamma distribution.
%We shall say that such a distribution is
\begin{defi}
\mylabel{defneargamma}
{\rm
 Let $\nu$ be a probability on  $\RR^+$  absolutely continuous with respect to the  Lebesgue
measure, with density $h$ and repartition function $$H(t)=\int_0^th(u)\;du\;.$$
Such a measure will be said to be \emph{nearly gamma} provided it satisfies the following set of conditions:
\bdes
\iti The set of $t \geq 0$ such that $h(t) > 0$ is an interval $I$;
\itii $h$ restricted to $I$ is continuous;
\itiii The map  $\psi : I \mapsto \RR$ defined by
$$\psi(y)=\frac{g\circ G^{-1}(H(y))}{h(y)}\;,$$
is such that :
\bdes
\ita There exists a positive real number $A$ such that
$$\forall y\in I,\;\psi(y)\leq
A\sqrt{y}\;.$$
\itb There exists $\eps >0$ such that, as $a$ goes to zero,
$$\nu\left( y\mbox{ :  }\psi(y)\leq a\right)=O(a^\eps)\;.$$
\edes
\edes
}
\end{defi}
In Definition \ref{defneargamma}, condition $(iii)$ is of course the most
tedious to check.  A simple sufficient condition for a
probability measure to be nearly gamma will be given in Lemma
\ref{lemneargamma}, the proof of which relies on the following asymptotics for
the Gaussian repartition function $G$.
\begin{lemm}
\label{lemasymptotgG}As $x$ tends to $-\infty$,
$$G(x)=g(x)\left(\frac{1}{|x|}+o\left(\frac{1}{x}\right)\right)\;,$$
and as $x$ tends to $+\infty$,
$$G(x)=1-g(x)\left(\frac{1}{x}+o\left(\frac{1}{x}\right)\right)\;.$$
Consequently,
$$g\circ G^{-1}(y)\stackrel{y\rightarrow 0}{\sim}y\sqrt{-2\log y}\;,$$
and
$$g\circ G^{-1}(y)\stackrel{y\rightarrow 1}{\sim}(1-y)\sqrt{-2\log (1-y)}\;.$$
\end{lemm}
\begin{dem}
A simple change of variable $u=x-t$ in $G$ gives:
$$G(x)=g(x)\int_{0}^{+\infty}e^{-\frac{t^2}{2}+xt}\;dt\;,$$
Integrating by parts, we get:
\begin{eqnarray*}
G(x)&=&g(x)\left(-\frac{1}{x}+\frac{1}{x}\int_{0}^{+\infty}te^{-\frac{t^2}{2}+xt}\;dt\right)\;,\\
&=&g(x)\left(\frac{1}{|x|}+o\left(\frac{1}{x}\right)\right)\;,\\
\end{eqnarray*}
as $x$ goes to $-\infty.$ Since
$G(-x)=1-G(x)$, we get that, as $x$ goes to $+\infty$:
$$G(x)=1-g(x)\left(\frac{1}{x}+o\left(\frac{1}{x}\right)\right)\;.$$
Let us turn to the asymptotic of $g\circ G^{-1}(y)$ as $y$ tends to zero. Let
$x=G^{-1}(y)$, so that ``$y$ tends to zero'' is equivalent to ``$x$ tends to
$-\infty$''. One has therefore,
\begin{eqnarray*}
G(x)&=&\frac{g(x)}{|x|}(1+o(1))\;,\\
\log G(x)&=&\log g(x)-\log |x| +o(1)\;,\\
&=&-\frac{x^2}{2}-\log |x| +O(1)\;,\\
\log G(x)&=&-\frac{x^2}{2}(1+o(1))\;,\\
|x|&=&\sqrt{-2\log G(x)}\;.\\
\end{eqnarray*}
Since $g(x)=|x|G(x)(1+o(1))$,
$$g(x)=G(x)\sqrt{-2\log G(x)}(1+o(1))\;,$$
and therefore,
$$g\circ G^{-1}(y)=y\sqrt{-2\log y}(1+o(1))\;,$$
as $y$ tends to zero. The asymptotic of $g\circ
G^{-1}(y)$ as $y$ tends to 1 is derived in the same way.
\end{dem}
Given two functions $r$ and $l$, we write $l(x) = \Theta(r(x))$ as $x$ goes to $x*$ provided there exist positive constants $C_1 \leq C_2$ such that
$$C_1 \leq \liminf_{x \rightarrow x*}  \frac{r(x)}{l(x)} \leq \limsup_{x \rightarrow x*} \frac{r(x)}{l(x)} \leq C_2.$$

\begin{lemm}
\label{lemneargamma}
Assume that condition $(i)$ and $(ii)$ of  Definition \ref{defneargamma} hold. Let $0 \leq \underline{\nu} < \overline{\nu} \leq \infty$ denote the endpoints of $I.$ Assume furthermore condition $(iii)$ is replaced by conditions $(iv)$ and $(v)$ below.
\bdes
\itiv There exists $\alpha >-1$ such that as $x$ goes to $\underline{\nu}$,
$$h(x)=\Theta\left((x-\underline{\nu})^\alpha\right)\;,$$
\itv $\overline{\nu} < \infty$ and there exists $\beta >-1$ such that as $x$ goes to $\overline{\nu}$,
$$h(x)=\Theta\left((\overline{\nu}-x)^\beta\right)\;,$$
\noindent or $\overline{\nu} = \infty$ and
$$\exists A>\underline{\nu},\;\forall t\geq A,\; C_1h(t)\leq \int_t^{\infty} h(u)\;du\leq C_2h(t)\;,$$
\noindent where $C_1$ and $C_2$ are positive constants.
\edes
Then, $\nu$ is nearly gamma.
\end{lemm}
\begin{dem}
Since $h$ is a continuous function on $]\underline{\nu},\overline{\nu}[$, it
attains its minimum on every compact set included in
$]\underline{\nu},\overline{\nu}[$. The minimum of $h$ on
$[a,b]$ is therefore strictly positive as
soon as $\underline{\nu}<a\leq b<\overline{\nu}$. In order to show
that condition $(iii)$ holds, we thus have to concentrate on the behaviour of the
function $\psi$ near $\underline{\nu}$ and $\overline{\nu}$. Condition $(iv)$
implies that, as $x$ goes to $\underline{\nu}$,
\begin{equation}
\label{DLHgauche}
H(x)=\Theta\left((x-\underline{\nu})^{\alpha+1}\right)\;.
\end{equation}
This, via Lemma \ref{lemasymptotgG}, leads to
\begin{equation}
\label{DLpsigauche}
\psi(x)=\Theta\left((x-\underline{\nu})\sqrt{-\log(x-\underline{\nu})}
  \right)\;,
\end{equation}
as $x$ goes to $\underline{\nu}$. Similarly, if $\overline{\nu} < \infty$,
  condition $(v)$ implies that, as $x$ goes to $\overline{\nu}$,:
\begin{equation}
\label{DLHdroite}
H(x)=\Theta\left((\overline{\nu}-x)^{\beta+1}\right)\;,
\end{equation}
which leads via Lemma \ref{lemasymptotgG} to
\begin{equation}
\label{DLpsidroite}
\psi(x)=\Theta\left((\overline{\nu}-x)\sqrt{-\log(\overline{\nu}-x)}
  \right)\;,
\end{equation}
as $x$ goes to $\overline{\nu}$. Therefore, if $\overline{\nu} < \infty$,
  part $(a)$ of condition $(iii)$ holds. Besides, asymptotics
  (\ref{DLpsigauche}) and (\ref{DLpsidroite}) imply that for $a$ small enough,
\begin{eqnarray*}
\nu\left( y\mbox{ s.t  }\psi(y)\leq
  a\right)&\leq&\PP((x-\underline{\nu})\sqrt{-\log(x-\underline{\nu})}\leq
Ca)\\
&&+
\PP((\overline{\nu}-x)\sqrt{-\log(\overline{\nu}-x)}\leq Ca)\;,
\end{eqnarray*}
for a positive constant $C$. Thus,
\begin{eqnarray*}
\nu\left( y\mbox{ :  }\psi(y)\leq a\right)& \leq & \PP((x-\underline{\nu})\leq
Ca)+\PP((\overline{\nu}-x)\leq
Ca)\;,\\
& =&H(\underline{\nu}+Ca)+1-H(\overline{\nu}-Ca)\;,\\
&=& O(a^{\alpha+1})+O(a^{\beta+1})\;,
\end{eqnarray*}
as $a$ goes to zero, according to equations (\ref{DLHgauche}) and
(\ref{DLHdroite}). Therefore, condition $(iii)$ holds when $\overline{\nu}$ is
finite.

Now, suppose that $\overline{\nu}=\infty$. Condition $(v)$ implies:
$$\forall t\geq A,\;1/C_2\leq \frac{h(t)}{\int_t^{\infty} h(u)}\;du\leq 1/C_1\;.$$
Integrating this inequality between A and $y$ leads to the existence of three
positive constants $B$, $C'_1$ and $C'_2$ such that:
$$\forall y\geq B,\;C'_1y\leq \log\frac{1}{1-H(y)}\leq C'2y\;.$$
Thus,
\begin{equation}
\label{DLpsiinfini}
\forall y\geq B,\;C_1\sqrt{C'_1y}\leq \psi(y)\leq C_2\sqrt{C'_2y}\;.
\end{equation}
This, combined with equation (\ref{DLpsigauche}) proves that part $(a)$ of
condition $(iii)$ holds. Now, equations (\ref{DLpsigauche}) and
(\ref{DLpsiinfini}) imply that for $a$ small enough, there is a constant $C$
such that,
$$\nu\left( y\mbox{ :  }\psi(y)\leq
  a\right)\leq\PP((x-\underline{\nu})\sqrt{-\log(x-\underline{\nu})}\leq
Ca)\;,$$
which was already proved to be of order $O(a^{\alpha+1})$. This concludes the proof of Lemma \ref{lemneargamma}.
\end{dem}
\brem
With the help of Lemma \ref{lemneargamma}, it is easy to check that most usual
distributions are nearly gamma. This includes all gamma and beta
distributions, as well as any probability measure whose
density is bounded away from 0 on its support, and notably the uniform
distribution on
$[a,b]$, with $0\leq a<b$. Nevertheless, remark that some distributions which 
have a sub-exponential upper tail may not satisfy the assumptions of Lemma
\ref{lemneargamma}, and be nearly gamma, though. For example, this is the case of the
distribution of $|N|$, where $N$ is a standard Gaussian random variable.
\erem
Now, we can state the main result of this article.
\begin{theo}
\label{theoFPPexpo} Let $\nu$ be a nearly gamma probability measure, with finite
moment of order 2. Let  $\mu$
denote the measure $\nu^{\otimes E}$. Then,
$$\Var_\mu(f_v)= O\left(\frac{|v|}{\log |v|}\right)\;,$$
as $|v|$ tends to infinity.
\end{theo}
In order to prove Theorem \ref{theoFPPexpo}, we use the same averaging
argument as in Benjamini et al. \cite{BenjaminiKalaiSchramm03}. It relies on the following lemma:
\begin{lemm}
\label{lemaveraging}
There exists a constant $c>0$, such that, for every $m\in \NN^*$, there exists a
function $g_m$ from $\{0,1\}^{m^2}$ to $\{0,\ldots,m\}$ such that:
$$\forall q\in\{1,\ldots ,m^2\},|\nabla_qg_m|\in\{0,1/2\}\;,$$
and
$$\max_{y\in\{0,\ldots,m\}}\lambda(x\mbox{ s.t. }g_m(x)=y)\leq\frac{c}{m}\;.$$
\end{lemm}
Since Benjamini et al. do not give a proof for this lemma, we
offer the following one.
\begin{dem}
From Stirling's Formula,
$$\binom{m^2}{\lfloor m^2/2\rfloor}.
\frac{m}{2^{m^2}}\xrightarrow[n\rightarrow \infty]{}\frac{1}{\sqrt{2\pi}}\;,$$
and this implies that the following supremum is finite:
$$c_1=\sup\left\lbrace 2\binom{m^2}{\lfloor
      m^2/2\rfloor}.\frac{m}{2^{m^2}}\mbox{ s.t. }m\in\NN^*\right\rbrace\;.$$
Notice also that $c_1\geq 1$. Now, let $\preceq$ denote the alphabetical order
      $\{0,1\}^{m^2}$, and let us list the elements in $\{0,1\}^{m^2}$ as
      follows:
$$(0,0,\ldots ,0)=x_1\preceq x_1\preceq \ldots \preceq x_{2^{m^2}}=(1,1,\ldots ,1)\;.$$
For any $m$ in $\NN^*$, we define the
      following integer:
$$k(m)=\left\lceil\frac{2^{m^2}}{m}\right\rceil\;,$$
and the following function on $\{0,1\}^{m^2}$:
$$\forall i\in \{1,\ldots
,2^{m^2}\},\;g_m(x_i)=\left\lfloor\frac{i}{k(m)}\right\rfloor\;.$$
Remark that $g_m(x_{2^{m^2}})\leq m/c_1\leq 1$. Therefore, $g$ is a function
from $\{0,1\}^{m^2}$ to $\{0,\ldots m\}$. Now, suppose that $x_i$ and $x_l$
differ from exactly one coordinate. Then, 
\begin{eqnarray*}
|i-l|&\leq & \binom{m^2}{i}+\binom{m^2}{l}\;,\\
&\leq & 2\binom{m^2}{\lfloor m^2/2\rfloor}\;,\\
&\leq & c_1\frac{2^{m^2}}{m}\;,\\
&\leq & k(m)\;.
\end{eqnarray*}
Consequently,
\begin{eqnarray*}
g_m(x_i)-g_m(x_l)&\leq&\left\lfloor\frac{l}{k(m)}+1\right\rfloor-\left\lfloor\frac{l}{k(m)}+1\right\rfloor\;,\\
&=& 1\;,
\end{eqnarray*}
which implies that $|\nabla_qg|\in\{0,1/2\}$. Finally, for any $y\in\{0,\ldots
m\}$, $g$ takes the value $y$ at most $k(m)$ times, and
\begin{eqnarray*}
\lambda(x\mbox{ s.t. }g_m(x)=y)& \leq & \frac{k(m)}{2^{m^2}}\;,\\
&\leq &\frac{c_1}{m}+\frac{1}{2^{m^2}}\;,\\
&\leq &\frac{2c_1}{m}\;.
\end{eqnarray*}
So the lemma holds with $c=2c_1$.
\end{dem}
We are now ready to prove Theorem \ref{theoFPPexpo}.

\vskip 2mm\noindent {\it Proof of Theorem \ref{theoFPPexpo}} :~
In the whole proof, $Y$ shall denote
a random variable with distribution $\nu$. We denote its second moment by $\sigma^2$:
$$\sigma^2=\EE(Y^2)\;.$$
 Let $m=\lceil
|v|^{\frac{1}{4}}\rceil$, and $S=\{1,\ldots,d\}\times\{1,\ldots,m^2\}$. Let
$c>0$ and $g_m$ be as in Lemma \ref{lemaveraging}. As in
\cite{BenjaminiKalaiSchramm03}, we define a random vertex in $\{0,\ldots ,m\}^d$ by the
following mean. For any $a=(a_{i,j})_{(i,j)\in S}\in \{0,1\}^S$, let
$$z=z(a)=\sum_{i=1}^dg_m(a_{i,1},\ldots ,a_{i,m^2})\mathbf{e_i}\;,$$
where $(\mathbf{e_1},\ldots ,\mathbf{e_d})$ denotes the standard basis of
$\ZZ^d$. We now equip the space $\{0,1\}^S\times \RR_+^E$ with the probability
measure $\lambda\otimes \mu$, where $\lambda$ is the uniform
measure on $\{0,1\}^S$, and we define the following function $\tilde{f}$ on $\{0,1\}^S\times
\RR_+^E$:
$$\forall (a,x)\in \{0,1\}^S\times \RR_+^E, \;
\tilde{f}(a,x)=d_x(z(a),v+z(a))\;.$$
The first important point to notice is that $f$ and $\tilde{f}$ are not too
far apart. Indeed, let $\alpha(a)$ be a path from 0 to $z(a)$,
such that $|\alpha(a)|=|z(a)|$ (here, $|\alpha|$ is the number of edges in
$\alpha$). Let $\beta(a)$ denote the
path $v+\alpha(a)$, which goes from $v$ to $v+z(a)$. Then: 
\begin{eqnarray*}
|\tilde{f}(a,x)-f(x)|&\leq& d_x(0,z(a))+d_x(v,v+z(a))\;,\\
&\leq &\sum_{e\in\alpha (a)}x_e+\sum_{e\in\beta (a)}x_e\;,\\
\end{eqnarray*}
Thus, using $|z|\leq m$,
$$ \|f - \tilde{f}\|_2 \leq  \left\|\sum_{e\in\alpha (a)}x_e\right\|_2+\left\|\sum_{e\in\beta (a)}x_e\right\|_2 \leq 2
\EE(|z(a)| \sigma) \leq 2m\sigma.$$
Therefore,
\begin{eqnarray}
%\Var(f)&=&\EE_{\lambda\otimes\mu}\left\lbrack\left(f-\tilde{f}
%  +\tilde{f}-\EE_{\lambda\otimes\mu}(\tilde{f})+\EE_{\lambda\otimes\mu}(\tilde%{f})-\EE_{\mu}(f)\right)^2\right\rbrack\;,\\
\nonumber \sqrt{\Var_\mu(f)} &=& \|f - \EE(f)\|_2 \leq \|f - \tilde{f}\|_2 + \|\tilde{f} - \EE(\tilde{f})\|_2 + \|\EE(\tilde{f})-\EE(f)\|_2 \;\\
\label{majovarfftilde} &\leq&  2 \|f - \tilde{f}\|_2 + \|\tilde{f} - \EE(\tilde{f})\|_2
\leq 4m\sigma + \sqrt{\Var(\tilde{f}) }
\end{eqnarray}
It remains to bound $\Var(\tilde{f})$. To this end, we will
use Corollary \ref{corochgtvar}. Using the notations of Corollary \ref{corochgtvar} and Definition
\ref{defneargamma}, we need to prove that $\tilde{f}\circ
(Id,\widetilde{H^{-1}\circ G})$ belongs to
$H_1^2(\lambda\otimes\gamma^\NN)$ when $\nu$ is nearly gamma. This, and the
application of Corollary \ref{corochgtvar} to $\tilde{f}$ is contained in the following Lemma.
\begin{lemm}
\label{lemderiveefv}
 If $\nu$ is nearly gamma, the function $f_v\circ\widetilde{H^{-1}\circ G}$ belongs to
  $H_1^2(\gamma^\NN)$, $\tilde{f}\circ
(Id,\widetilde{H^{-1}\circ G})$ belongs to
  $H_1^2(\lambda\otimes\gamma^\NN)$ and,
$$Var(\tilde{f})\leq\sum_{q\in S}\left\|\nabla_q \tilde{f}\right\|_2^2+ 2
\sum_{e\in E}\left\|\nabla_e\tilde{f}\right\|_2^2\phi\left(\frac{\left\|\nabla_e
      \tilde{f}\right\|_1}{\left\|\nabla_e \tilde{f}\right\|_2}\right)\;,
$$
where for every edge $e$,
\begin{equation}
\label{eqapplicoroFPP}
\nabla_e\tilde{f}(a,x)=\psi(x_e)\frac{\partial
  \tilde{f}}{\partial x_e}(a,x)\;.
\end{equation}
Furthermore, conditionally to $z$, there is almost surely only one $x$-geodesic from $z$ to
$z+v$, denoted by $\gamma_x(z)$, and:
$$\frac{\partial
  \tilde{f}}{\partial x_e}(a,x)=\II_{e\in\gamma_x(z(a))}\;.$$
\end{lemm}
\begin{dem}
The fact that $f_v\circ\widetilde{H^{-1}\circ G}$ and $\tilde{f}\circ
(Id,\widetilde{H^{-1}\circ G})$ are in $L^2$ follows from the basic fact that
$\nu$ has a finite moment of order 2, and that $d_x(u,v)$ is dominated by a sum of $|v-u|$ independent variables with
distribution $\nu$, which is the length of a deterministic path of length $|v-u|$. We shall prove that $f_v\circ\widetilde{H^{-1}\circ G}$ satisfies the integration by part formula {\bf
  (a)} of the definition of $H_1^2$. The similar result for $\tilde{f}$ is
obtained in the same way.

Let $e$ be an edge of $\ZZ^d$, and $x^{-e}$ be an element of
$(\RR^+)^{E(\ZZ^d)\setminus\{e\}}$. For every $y$ in $\RR^+$, we denote by
$(x^{-e},y)$ the element $x$ of $(\RR^+)^{E(\ZZ^d)}$ such that:
$$x_e=y\mbox{ and }\forall e'\not = e,\;x_{e'}=x^{-e}_{e'}\;.$$
Now, we fix $x^{-e}$ in $(\RR^+)^{E(\ZZ^d)}$. We denote by $g_e$ the function
defined on $\RR^+$ by:
$$g(y)=f_v(x^{-e},y)\;.$$
We will show that there is a nonnegative real number $y_\infty$ such that:
\begin{equation}
\label{claimg}
\left\lbrace\begin{array}{l}\forall y\leq y_\infty,\;g(y)=g(0)+y\\
\mbox{ and }\\
\forall y > y_\infty,\;g(y)=g(y_\infty)\end{array}\right.
\end{equation}
For any $n\geq |v|$, let us denote by $\Gamma_n$ the set of paths from 0 to $v$ whose number of
edges is not greater than $n$. We have:
$$g(y)=\inf_{n\geq |v|}g_n(y)\;,$$
where 
$$g_n(y)=\inf_{\gamma\in\Gamma_n}\sum_{e'\in \gamma}(x^{-e},y)_{e'}\;.$$
The functions $g_n$ form a nonincreasing sequence of nondecreasing functions:
$$\forall n\geq|v|,\;\forall y\in\RR^+,\;\forall y'\geq y,\;g_{n+1}(y)\leq
g_n(y)\leq g_n(y')\;.$$
In particular, this implies that for every $y$ in $\RR^+$,
$$g(y)=\lim_{n\rightarrow \infty}g_n(y)\;.$$ 
Now, we claim that, for every $n\geq |v|+3$, there exists $y_n\in\RR^+$ such that:
\begin{equation}
\label{claimgn1}
\left\lbrace\begin{array}{l}\forall y\leq y_n,\;g_n(y)=g_n(0)+y\\
\mbox{ and } \\
\forall y > y_n,\;g_n(y)=g_n(y_n)\end{array}\right.
\end{equation}
and furthermore, 
\begin{equation}
\label{claimgn2}
\mbox{ the sequence }(y_n)_{n\geq|v|+3}\mbox{ is nonincreasing.}
\end{equation}
Indeed, since
$\Gamma_n$ is a finite set, the infimum in the definition of $g_n$
is attained. Let us call a path which attains this infimum an
$(n,y)$-geodesic and let $\tilde{\Gamma}(n,y,e)$ be the set of $(n,y)$-geodesics
which contain the edge $e$. Remark that as soon
as $n\geq |v|+3$, there
exists a real number $A$ such that $e$ does not belong to any $(n,A)$-geodesic:
it is enough to take $A$ greater than the sum of the length of three edges forming a path between the
end-points of the edge $e$. Therefore, the following supremum is finite:
$$y_n=\sup\{y\in\RR^+\mbox{ s.t. }\tilde{\Gamma}(n,y,e)\not =\emptyset\}\;.$$
Now, if $e$ belongs to an $(n,y)$-geodesic $\gamma$, for any
$y'\leq y$, $\gamma$ is an $(n,y')$-geodesic to which $e$ belongs, and $g_n(y)-g_n(y')=y-y'$. If
$\tilde{\Gamma}(n,y,e)$ is empty, then for any $y'\geq y$, $e$ does
not belong to any $(n,y')$-geodesic, and $g_n(y)=g_n(y')$. This proves that:
$$\forall y< y_n,\;g_n(y)=g_n(0)+y\;,$$
$$\forall y,y' > y_n,\;g_n(y)=g_n(y')\;.$$
Since $g_n$ is continuous, we have proved claim (\ref{claimgn1}). Now remark
that if $e$ does not belong to any $(n,y)$-geodesic, then $e$ does not belong to
any $(n+1,y)$-geodesic, since $\Gamma_n\subset\Gamma_{n+1}$. Therefore,
$y_{n+1}\leq y_n$, and this proves claim (\ref{claimgn2}). Since
$(y_n)_{n\geq|v|+3}$ is nonnegative, it converges to a nonnegative number
$y_\infty$ as $n$ tends to infinity. Now, let $n$ be a integer greater than
$|v|+3$:
$$\forall n\geq N,\;\forall y,y'>y_n,\;g_n(y)=g_n(y')\;.$$
Since $y_n\leq y_N$, 
$$\forall n\geq N,\;\forall y,y'>y_N,\;g_n(y)=g_n(y')\;.$$
Letting $n$ tend to infinity in the last equation, we get:
$$\forall N\geq |v|+3,\;\forall y,y'>y_N,\;g(y)=g(y')\;.$$
Therefore,
$$\forall y,y'>y_\infty,\;g(y)=g(y')\;.$$
On the other side,
$$\forall n\geq |v|+3,\;\forall y\leq y_n,\;g_n(y)=g_n(0)+y\;.$$
Since $y_n\geq y_\infty$,
$$\forall n\geq |v|+3,\;\forall y\leq y_\infty,\;g_n(y)=g_n(0)+y\;.$$
Letting $n$ tend to infinity in the last expression, we get:
$$\forall y\leq y_\infty,\;g(y)=g(0)+y\;.$$
Finally, $g$ is continuous. Indeed, the convergent sequence $(g_n)$ is
uniformly equicontinuous, since all these functions are
1-Lipschitz. Furthermore, they are uniformly bounded, since they are
nonnegative, and less than $g_{|v|+3}(y_{|v|+3})$. Then, the continuity of
$g$ follows from Arzelà-Ascoli Theorem. We have proved claim
(\ref{claimg}). Remark that $y_\infty=y_\infty(x^{-e})$ depends on
$x^{-e}$. We define, for any $x^{-e}$,
$$h_e(x^{-e},x_e)=\left\lbrace\begin{array}{l}1\mbox{ if } x_e\leq y_\infty(x^{-e}) \\ 0\mbox{ if } x_e >y_\infty(x^{-e})\end{array}\right.\;.$$
It is easy to see that, for any
smooth function $F:\RR\rightarrow \RR$ having compact support, for any
$x^{-e}$,
\begin{equation}
\label{eqIPP}
-\int_{\RR} F'(x_e)f_v(x^{-e},x_e)\;dx_e=\int_{\RR}
F(x_e)h_e(x^{-e},x_e)\;dx_e\;.
\end{equation}
It is known that there is almost surely a geodesic from 0 to $v$ (see
\cite{Howard04} for instance), i.e the
infimum in the definition of $f_v$ is attained with probability
1. Furthermore, in this setting, where the distribution of the lengths is
continuous, there is almost surely only one $x$-geodesic from $0$ to
$v$. For any  $x$, we shall denote by $\gamma_x(0)$ the unique $x$-geodesics from
  $0$ to $0+v$. Then, with $\nu$-probability 1, one can see from the definitions
  of $y_n$ and $y_\infty$ that:
\begin{equation}
\label{eqexpressionder}h_e(x^{-e},x_e)=\II_{e\in\gamma_x(0)}\;.
\end{equation}
Performing the change of variable $x\mapsto \widetilde{H^{-1}\circ G}$ in
equation (\ref{eqIPP}), one gets the integration by parts formula {\bf (a)}
for $f_v\circ \widetilde{H^{-1}\circ G}$, with the following partial derivative with
respect to $x_e$:
$$x\mapsto \psi(x_e)h_e(\widetilde{H^{-1}\circ G}(x))\;.$$
The fact that the sum of the $L^2$ norms of the partial derivatives of $f_v$
and $\tilde{f}$ is
finite follows from estimate (\ref{majsomnorm}). Then, inequality
(\ref{eqapplicoroFPP}) is a consequence of Corollary
\ref{corochgtvar}, and the expression of $\frac{\partial
  \tilde{f}}{\partial x_e}(a,x)$ is derived in the same way than (\ref{eqexpressionder}).
\end{dem}
Our next goal is to find a good upper-bound for $\left\|\nabla_e
  \tilde{f}\right\|_{2}^2$, for any edge $e$. Let $A$ be as
  in Definition \ref{defneargamma}. It follows from Lemma \ref{lemderiveefv}
  and the fact that $\nu$ is nearly gamma that:
\begin{eqnarray*}
\left\|\nabla_e
  \tilde{f}\right\|_{2}^2&=&\EE\left(\psi(x_e)^2\II_{e\in\gamma_x(z)}\right)\;,\\
&\leq&A^2\EE_{\lambda}\left(\EE_{\mu}\left(x_e\II_{e\in\gamma_x(z)}\right)\right)\;.
\end{eqnarray*}
Now, we use the fact that for any fixed $z$, $\mu$ is invariant under translation by $z$.
\begin{eqnarray*}
\left\|\nabla_e
  \tilde{f}\right\|_{2}^2&\leq&
%A^2\EE_{\lambda}\left(\EE_{\mu}\left(x_{e-z}\II_{e\in\gamma_{x+z}(z)}\right)\right)\;,\\
%&=&
A^2\EE_{\lambda}\left(\EE_{\mu}\left(x_{e-z}\II_{e-z\in\gamma_{x}(0)}\right)\right)\;,\\
&=&A^2\EE_{\mu}\left(\sum_{e'\in\gamma_{x}(0)}\EE_{\lambda}\left(x_{e'}\II_{e-z=e'}\right)\right)\;,\\
&=&A^2\EE_{\mu}\left(\sum_{e'\in\gamma_{x}(0)}x_{e'}\PP_{\lambda}(z=e-e')\right)\;,\\
&\leq& A^2\sup_{z_0}\PP(z=z_0)\EE_{\mu}\left(\sum_{e'\in\gamma_{x}(0)\cap \mathcal{Q}_e}x_{e'}\right)\;,\\
\end{eqnarray*}
where $\mathcal{Q}_e=\left\lbrace e'\in E\left(\ZZ^d\right)\mbox{
    s.t. }\PP(z=e-e')>0\right\rbrace $. Using Lemma \ref{lemaveraging}, 
$$ \sup_{z_0}\PP(z=z_0)\leq \left(\frac{c}{m}\right)^d\;.$$
Let $e^-$ and $e^+$ denote the end-points of $e$, and $\mathcal{B}(0,dm)$ be the $L^1$-ball with center
0 and radius $dm$. Remark now that $\mathcal{Q}_e$ is included in the set of edges
$e+\mathcal{B}(0,dm)$, which is itself included in $\mathcal{B}(e^-,dm+1)$.  This simply follows from the fact that $g_m$ takes its
values in $\{0,\ldots , m\}$. Let $r$ be a deterministic path going through
every vertex of the surface of the ball $\mathcal{B}(e^-,dm+1)$, and such that
there is a constant $C$ (depending only on $d$) such that $|r|\leq Cm^{d-1}$. From the
definition of a geodesic, we get:

 %% Let $v_1$ be the
%%     first vertex of $\gamma$ belonging to $\mathcal{Q}_e$, when one goes from
%%     0 to $v$, and $v_2$ be the last point of $\gamma$ belonging to
%%     $\mathcal{Q}_e$. Let $\tilde{\gamma}$ be the restriction of $\gamma$ from
%%     $v_1$ to $v_2$. We have:

%% \begin{eqnarray*}
%% \sum_{e'\in\gamma\cap\mathcal{Q}_e}x_{e'}&=&\sum_{e'\in\tilde{\gamma}\cap\mathcal{Q}_e}x_{e'}\;,\\
%% &\leq&\inf_{\alpha:(v_1,v_2)}\sum_{e'\in\alpha}x_{e'}\;,\\
%% &\leq&\inf_{\substack{\beta:(v_1,v_2)\\ \beta\subset conv(\mathcal{Q}_e)}}\sum_{e'\in\alpha}x_{e'}\;,\\
%% \end{eqnarray*}
%% where $conv(\mathcal{Q}_e)$ is the intersection of the convex hull of
%% $\mathcal{Q}_e$ in $\RR^d$ with $\ZZ^d$. This leads to:

\begin{eqnarray*}
\EE_{\mu}\left(\sum_{e'\in\gamma_{x}(0)\cap
    \mathcal{Q}_e}x_{e'}\right)&\leq&\EE\left(\sum_{e'\in \mathcal{B}(e^-,dm+1)}x_{e'}\right)\;,\\
&\leq &\EE\left(\sum_{e'\in
   r}x_{e'}\right)\;,\\
&\leq &Cm^{d-1}\EE(Y)\;.
\end{eqnarray*}
Therefore, 
\begin{equation}
\label{majnorm2}
\left\|\nabla_e
  \tilde{f}\right\|_{2,\lambda\otimes\mu}^2 \leq C
  A^2dm^{d-1}\EE(Y)\left(\frac{c}{m}\right)^d \leq A^2d\EE(Y)\frac{c^d}{|v|^{\frac{1}{4}}}\;.
\end{equation}
In order to use Corollary \ref{coroTalFPPexpo}, one needs to bound from below the
  quotient $\left\|\nabla_e
  \tilde{f}\right\|_{2}/\left\|\nabla_e
  \tilde{f}\right\|_{1}$. Let $a$ be a positive real number. 
\begin{eqnarray*}
\left\|\nabla_e \tilde{f}\right\|_{1}&= &\EE\left(
\psi(x_e)\II_{e\in\gamma_x(z)}\right)\;,\\
&=&\EE\left(
\psi(x_e)\II_{e\in\gamma_x(z)}\II_{\psi(x_e)\leq
  a}\right)+\EE\left(
\psi(x_e)\II_{e\in\gamma_x(z)}\II_{\psi(x_e) >
  a}\right)\;,\\
&\leq&\EE\left(\psi(x_e)\II_{\psi(x_e)\leq
  a}\right)+\frac{1}{a}\EE\left(
\psi(x_e)^2\II_{e\in\gamma_x(z)}\II_{\psi(x_e) >
  a}\right)\;,\\
&\leq& a\nu(\psi(x_e)\leq a)+\frac{1}{a}\left\|\nabla_e \tilde{f}\right\|_{2}^2\;.
\end{eqnarray*}
Since $\nu$ is nearly gamma, there exists $\eps >0$ such that:
$$\nu(\psi(x_e)\leq a)=O(a^\eps)\;,$$
as $a$ goes to zero. Now, let us choose $a=\left\|\nabla_e
  \tilde{f}\right\|_{2}^{2/(2+\eps)}$.
Thus,
$$\left\|\nabla_e \tilde{f}\right\|_{1}=O\left(\left\|\nabla_e
  \tilde{f}\right\|_{2}^{1+\frac{\eps}{2+\eps}}\right)\;.$$
This, via inequality (\ref{majnorm2}), leads to
\begin{equation}
\label{majquotientnorm}
\log \frac{\left\|\nabla_e
  \tilde{f}\right\|_{1}}{\left\|\nabla_e
  \tilde{f}\right\|_{2}}=O(\log |v|)
\end{equation}
In addition, if one chooses a particular $x$-geodesic $\gamma$ from $z$ to
$z+v$, and a $L^1$-geodesic $\alpha(z)$ from $z$ to
$z+v$,
\begin{eqnarray}
\nonumber\sum_{e\in E}\left\|\nabla_e
  \tilde{f}\right\|_{2}^2& =
  &\EE\left(\sum_{e\in\gamma}\psi(x_e)^2 \right)\;,\\
\nonumber& \leq
  & A\EE\left(\sum_{e\in\gamma}x_e \right)\;,\\
\nonumber&\leq&A\EE(Y)\times\EE_{\lambda}\left(\left|\alpha(z)\right|\right)\;,\\
\label{majsomnorm}\sum_{e\in E}\left\|\nabla_e
  \tilde{f}\right\|_{2}^2&\leq&A\EE(Y)|v|\;.
\end{eqnarray}
Collecting estimates (\ref{majquotientnorm}) and (\ref{majsomnorm}),
\begin{eqnarray}
\nonumber\sum_{e\in E}\left\|\nabla_e\tilde{f}\right\|_2^2\phi\left(\frac{2\sqrt{2}\left\|\nabla_e
      \tilde{f}\right\|_1}{\pi\left\|\nabla_e\tilde{f}\right\|_2}\right)&\leq &\sum_{e\in E}\left\|\nabla_e\tilde{f}\right\|_2^2\sup_{e\in E}\phi\left(\frac{2\sqrt{2}\left\|\nabla_e
      \tilde{f}\right\|_1}{\pi\left\|\nabla_e\tilde{f}\right\|_2}\right)\;,\\
\label{majsomedges}&= &O\left(\frac{|v|}{\log |v|}\right)\;,
\end{eqnarray}
where we used the fact that $\phi(u) \sim -1/(2 \log(u))$ when $u$ goes to zero. According to Lemma \ref{lemaveraging}, for any $q\in S$,
$|\nabla_qg_m|\in\{0,1/2\}$. Therefore, 
$$\forall q\in S,\;\left\|\nabla_qf\right\|_2^2\leq\sigma^2\;,$$
and thus,
$$\sum_{q\in S}\left\|\nabla_qf\right\|_2^2 \leq |S|\sigma^2=dm^2\sigma^2\;.$$
Inequality (\ref{majsomedges}), the assumption $|m|=\lceil
|v|^{1/4}\rceil$ and Lemma \ref{lemderiveefv} lead therefore to:
$$\Var(\tilde{f})=O\left(\frac{|v|}{\log |v|}\right)\;.$$
This, together with inequality  (\ref{majovarfftilde}) implies that:
$$\Var(f_v)= O\left(\frac{|v|}{\log |v|}\right)\;.$$
The proof of Theorem \ref{theoFPPexpo} is complete.
 \hfill $\square$ \vskip 2mm \noindent

%%%%%%%%%%%%%%%%%%%%%%%%%%%%%%%%%%%%%%%%%%%%%%%%%

%% \bibliography{../zero_un}
%% \bibliographystyle{apalike}

\end{document}